%% file: pde.new.ref.tex
\documentclass{proc-l}

\usepackage{graphicx} 
\usepackage{amsmath}

\graphicspath{ {./images/} }

\newtheorem{theorem}{Theorem}[section]
\newtheorem{lemma}[theorem]{Lemma}

\theoremstyle{definition}

\theoremstyle{remark}

\numberwithin{equation}{section}

\begin{document}
\title{The Exact Solutions of Certain Linear Partial Difference Equations}

\author{Chun-Kai Hwang}
\address{Department of Computer Science and Information Engineering, National Taiwan University, Taipei, Taiwan, 10617, R.O.C.}
\email{chunkai.hwang@gmail.com}
\author{Tzon-Tzer Lu}
\address{Department of Applied Mathematics, National Sun Yat-Sen University, Kaohsiung, Taiwan, 80424, R.O.C.}
\email{ttlu@math.nsysu.edu.tw}


\date{February 18, 2025 }




\begin{abstract}
Difference equations have many applications and play an important role in numerical analysis, probability, statistics, combinatorics, computer science, quantum consciousness, etc. We first prove that the partial differential equation is equivalent to partial difference equation with an example of heat equation. Additionally, we use generating functions to find the exact solutions of some simple linear partial difference equations. Then we extend it to more general partial difference equations of higher dimensions and obtain their solutions. Notice that Theorem 4.2 could provide a mathematical framework for understanding how information within a black hole is encoded on its event horizon, a key concept in the black hole information paradox. In addition, we extend it to n-dimensional case, Theorem 4.4, the high-order partial difference equations (HOPDE). We conclude that using multivariable power series as generating function is a very efficient method to solve partial difference equations.
\end{abstract}

\maketitle

\section{Introduction}

Many phenomena or problems occur in the natural world can be formulated by
the mathematical model of partial difference equations, such as discrete
potential theories, random walk problems, finite difference schemes,
gambling problems, chemical structures, discrete heat diffusion
equations,...,etc. Additionally, partial difference equations are related to quantum consciousness. The arrangement of microtublin follows the Fibonacci series from the quantum Hopfield networks\cite{behrman2006microtubules},\cite{srivastava2016modelling}, 
and Orch. OR (Orchestrated Objective Reduction) theory
\cite{hameroff2014consciousness},
\cite{penrose2011consciousness}
proposed by Penrose and Hameroff,. Roughly speaking, a partial difference equation is a
functional relation involving functions with multivariate discrete variables
over a discrete set.
 \par
Next, we take random work problems and heat diffusion as examples.
Suppose that an ant creeps along the one dimensional straight line, x-axis.
Assume the probability it doesn't move in a given unit of time is d, the
probability that it moves one unit to the positive and negative directions
along the x-axis in a unit of time are p and q respectively. Let
$a(i,j)$ be the probability that the ant is at the point $x=i$ after $j$-th units
of time. Then solving $a(i,j)$ is equivalent to solving the following
difference equation 
\[
a(i,j+1)=pa(i-1,j)+da(i,j)+qa(i+1,j).
\]
If the ant starts from the origin at time 0, then we have initial conditions  
\[
  a(i,0) =
    \begin{cases}
      1, & \text{if i = 0,} \\
      0, & \text{otherwise.}
    \end{cases}       
\]

Next, we discuss the discrete version of the heat equation 
\[
\frac{\partial u}{\partial t}=b\frac{\partial ^{2}u}{\partial x^{2}}%
,t>0,-\infty <x<\infty .
\]
Consider the distribution of heat through a long rod. Assume that the rod is
long enough so that we can sample it with the discrete integer set. Let $U_{i,j}$
be the temperature at the position $i$ after $j$-th units of time. At
time $j$, if the temperature $U_{i-1,j}$ is higher than $U_{i,j}$, heat will
flow from the position $i-1$ to $i$. From the experiment, the increase is
proportional to the \begin{math}U_{i-1,j}-U_{i,j}\end{math}, say \begin{math}r(U_{i-1,j}-U_{i,j})\end{math},
where $r$ is a physical parameter depending on the properties of the rod.
Similarly, heat will flow from the point $i+1$ to $i$. Hence, the total effects
can be written as the following descrete Newton law of cooling 
\[
U_{i,.j+1}-U_{i,.j}=r(U_{i-1,.j}-U_{i,.j})+r(U_{i+1,.j}-U_{i,.j})=r\triangle
_{1}^{2}U_{i-1,j}.
\]
We can regard the equation above as a special case of the difference
equation 
\[
U_{i,j+1}=aU_{i-1,j}+bU_{i,j}+cU_{i+1,j},
\]
where i$\in Z,$ j$\in N\cup \{0\}.$
Although they have different interpretation, they correspond to the same type of equation.The exact solution of this equation will be computed in Theorem 2.1.

There are many methods developed to solve the difference equation, such as the symbolic method\cite{cheng1999general}, 
Lagrange method\cite{mickens2015difference},
generating functions\cite{berman2014introduction},
symbolic calculus\cite{cheng1999difference},
characteristic equations, matrices iterations\cite{cheng2003partial},
...,etc. We will use multivariable power series as generating function to solve the exact solutions of certain partial difference equations. We demonstrate its usage for some simple cases in section 2, then extend it to general case in section 3 and multidimensional cases in section 4. It can be seen that this method is very efficient to solve the partial difference equations.

\section{Some Simple Cases}

Consider the following linear partial difference equations:

\begin{equation}
U_{i,j+1}=aU_{i-1,j}+bU_{i,j}+cU_{i+1,j}  \label{eq1}
\end{equation}

\begin{equation}
U_{i+m,j+1}=c_{1}U_{i,j}+c_{2}U_{i+1,j}+...+c_{n}U_{i+(n-1),j}  \label{eq2}
\end{equation}

\begin{equation}
U_{i+1,j+1}=aU_{i,j+1}+bU_{i+1,j}+cU_{i,j}\text{ }  \label{eq3}
\end{equation}

We will obtain the explicit solution of $U_{i,j}$ under the assumption of
$U_{i,0}$ given for all $i\in Z$. That is, if we know all of the physical
quantity on the x-axis, we can predict the other quantity $%
U_{i,j} $ at any position in the upper plane.

\par
We first consider the equation(\ref{eq1}). Lin\cite{lin1997}
deals with more general nonhomogeneous case of this equation.

\begin{theorem}
Let

\[
U_{i,j+1}=aU_{i-1,j}+bU_{i,j}+cU_{i+1,j} 
\]
\end{theorem}

with i$\in Z,j\in N\cup \{0\}$.
\ Given $U_{i,0}=\Psi (i)$, then its solution is of the form

\[
U_{i,j}=\sum_{m=0}^{j}\sum_{n=0}^{m}\binom{j}{m}\binom{m}{n}%
a^{n}b^{m-n}c^{j-n}\Psi (i+j-m-n). 
\]

\proof
Let 
\[
h(x,y)=\sum_{i=-\infty }^{\infty }\sum_{j=0}^{\infty }U_{i,j}x^{i}y^{j}. 
\]
Then we have

\[
\sum_{i=-\infty }^{\infty }\sum_{j=0}^{\infty
}U_{i,j+1}x^{i}y^{j}=a\sum_{i=-\infty }^{\infty }\sum_{j=0}^{\infty
}U_{i-1,j}x^{i}y^{j}+b\sum_{i=-\infty }^{\infty }\sum_{j=0}^{\infty
}U_{i,j}x^{i}y^{j}+c\sum_{i=-\infty }^{\infty }\sum_{j=0}^{\infty
}U_{i+1,j}x^{i}y^{j}, 
\]

\begin{eqnarray*}
\frac{1}{y}[h(x,y)-\sum_{i=-\infty }^{\infty }\Psi (i)x^{i}]
&=&axh(x,y)+bh(x,y)+\frac{c}{x}h(x,y), \\
h(x,y)[1-(ax+b+\frac{c}{x})y] &=&\sum_{i=-\infty }^{\infty }\Psi (i)x^{i},
\end{eqnarray*}

\begin{eqnarray*}
h(x,y) &=&\frac{\sum_{i=-\infty }^{\infty }\Psi (i)x^{i}}{1-(ax+b+cx^{-1})y}
\\
&=&\sum_{i=-\infty }^{\infty }\Psi (i)x^{i}\sum_{j=0}^{\infty
}(ax+b+cx^{-1})^{j}y^{j} \\
&=&\sum_{i=-\infty }^{\infty }\Psi (i)x^{i}\sum_{j=0}^{\infty }[\sum_{\substack{
0\leq \alpha ,\beta ,\gamma \leq j, \\ \alpha +\beta +\gamma =j}} 
\binom{j}{\alpha \beta \gamma }a^{\alpha }b^{\beta }c^{\gamma }x^{\alpha
-\gamma }]y^{j} \\
&=&\sum_{j=0}^{\infty }\sum_{i=-\infty }^{\infty }\Psi (i)x^{i}[\sum_{\substack{
0\leq \alpha ,\beta ,\gamma \leq j, \\ \alpha +\beta +\gamma =j}} 
\binom{j}{\alpha \beta \gamma }a^{\alpha }b^{\beta }c^{\gamma }x^{\alpha
-\gamma }]y^{j}.
\end{eqnarray*}
Let $\alpha =n,\beta =m-n,\gamma =j-m$, then we have

\begin{eqnarray*}
h(x,y) &=&\sum_{j=0}^{\infty }\{\sum_{i=-\infty }^{\infty }\Psi
(i)x^{i}\sum_{m=0}^{j}\sum_{n=0}^{m}\binom{j}{n,m-n,j-n}%
a^{n}b^{m-n}c^{j-n}x^{n-j+m}\}y^{j} \\
&=&\sum_{j=0}^{\infty }\{\sum_{i=-\infty }^{\infty
}\sum_{m=0}^{j}\sum_{n=0}^{m}\binom{j}{n,m-n,j-n}a^{n}b^{m-n}c^{j-n}\Psi
(i)x^{i-j+m+n}\}y^{j}.
\end{eqnarray*}
Change of variable $k=i-j+m+n$ to obtain

\[
h(x,y)=\sum_{j=0}^{\infty }\{\sum_{i=-\infty }^{\infty
}\sum_{m=0}^{j}\sum_{n=0}^{m}\binom{j}{n,m-n,j-n}a^{n}b^{m-n}c^{j-n}\Psi
(k+j-m-n)x^{k}\}y^{j}. 
\]
Thus the solution is of the form 
\[
U_{i,j}=\sum_{m=0}^{j}\sum_{n=0}^{m}\binom{j}{m}\binom{m}{n}%
a^{n}b^{m-n}c^{j-n}\Psi (i+j-m-n). 
\]

\endproof%

\begin{theorem}
Let 
\[
U_{i+m,j+1}=c_{1}U_{i,j}+c_{2}U_{i+1,j}+...+c_{n}U_{i+(n-1),j} 
\]
with i$\in Z,j\in N\cup \{0\}$.
Given $U_{i,0}=\Psi (i)$, then its solution is of the form
\end{theorem}

\[
U_{i,j}=\sum_{s_{r\in I}}\binom{j}{s_{1}s_{2...}s_{n}}\Pi
_{r=1}^{n}c_{r}^{s_{r}}\Psi (i+\sum_{r=1}^{n}rs_{r}-(m+1)j), 
\]

where $I=\{s_{1},s_{2},...,s_{n}|0\leq s_{k}\leq j,\sum_{k=1}^{n}s_{k}=j\}.$

\proof%
Let 
\[
h(x,y)=\sum_{i=-\infty }^{\infty }\sum_{j=0}^{\infty }U_{i,j}x^{i}y^{j}. 
\]
Then we have

\begin{eqnarray*}
&&\sum_{i=-\infty }^{\infty }\sum_{j=0}^{\infty }U_{i+m,j+1}x^{i}y^{j} \\
&=&c_{1}\sum_{i=-\infty }^{\infty }\sum_{j=0}^{\infty
}U_{i,j}x^{i}y^{j}+c_{2}\sum_{i=-\infty }^{\infty }\sum_{j=0}^{\infty
}U_{i+1,j}x^{i}y^{j}+...+c_{n}\sum_{i=-\infty }^{\infty }\sum_{j=0}^{\infty
}U_{i+(n-1),j}x^{i}y^{j},
\end{eqnarray*}

\[
\frac{1}{x^{m}y}[h(x,y)-\sum_{i=-\infty }^{\infty }\Psi (i)x^{i}%
]=\sum_{t=1}^{n}c_{t}x^{1-t}h(x,y), 
\]

\[
h(x,y)(1-\sum_{t=1}^{n}c_{k}x^{m-t+1}y)=\sum_{i=-\infty }^{\infty }\Psi
(i)x^{i}, 
\]

\begin{eqnarray*}
h(x,y) &=&\frac{\sum_{i=-\infty }^{\infty }\Psi (i)x^{i}}{%
1-(c_{1}x^{m}+c_{2}x^{m-1}+...+c_{n}x^{m-(n-1)})y} \\
&=&\sum_{i=-\infty }^{\infty }\Psi (i)x^{i}\sum_{j=0}^{\infty
}(c_{1}x^{m}+c_{2}x^{m-1}+...+c_{n}x^{m-(n-1)})^{j}y^{j} \\
&=&\sum_{i=-\infty }^{\infty }\Psi (i)x^{i}\sum_{j=0}^{\infty }\sum_{s_{r\in
I}}\binom{j}{s_{1}s_{2...}s_{n}}%
(c_{1}^{s_{1}}x^{ms_{1}})(c_{2}^{s_{2}}x^{(m-1)s_{2}})...(c_{n}^{s_{n}}x^{[m-(n-1)]s_{n}})y^{j}
\\
&=&\sum_{i=-\infty }^{\infty }\Psi (i)x^{i}\sum_{j=0}^{\infty }\sum_{s_{r\in
I}}\binom{j}{s_{1}s_{2...}s_{n}}\Pi _{r=1}^{n}c_{r}^{s_{r}}x^{a}y^{j} \\
&=&\sum_{j=0}^{\infty }\{\sum_{i=-\infty }^{\infty }\Psi
(i)x^{i}\sum_{s_{r\in I}}\binom{j}{s_{1}s_{2...}s_{n}}\Pi
_{r=1}^{n}c_{r}^{s_{r}}x^{a}\}y^{j} \\
&=&\sum_{j=0}^{\infty }\{\sum_{i=-\infty }^{\infty }\sum_{s_{r\in I}}\binom{j%
}{s_{1}s_{2...}s_{n}}\Pi _{r=1}^{n}c_{r}^{s_{r}}\Psi (i)x^{a+i}\}y^{j},
\end{eqnarray*}
where $a=\sum_{r=1}^{n}[m-(r-1)]s_{r}$, $I=\{s_{1},s_{2},...,s_{n}|0\leq
s_{k}\leq j$, $\sum_{k=1}^{n}s_{k}=j\}$. Let 
\[
k=a+i. 
\]
Then we have

\begin{eqnarray*}
h(x,y) &=&\sum_{j=0}^{\infty }\sum_{k=-\infty }^{\infty }\sum_{s_{r\in I}}%
\binom{j}{s_{1}s_{2...}s_{n}}\Pi _{r=1}^{n}c_{r}^{s_{r}}\Psi (k-a)x^{k}y^{j}
\\
&=&\sum_{j=0}^{\infty }\sum_{k=-\infty }^{\infty }\sum_{s_{r\in I}}\binom{j}{%
s_{1}s_{2...}s_{n}}\Pi _{r=1}^{n}c_{r}^{s_{r}}\Psi
(k-\sum_{r=1}^{n}[m-(r-1)]s_{r})x^{k}y^{j} \\
&=&\sum_{j=0}^{\infty }\sum_{k=-\infty }^{\infty }\sum_{s_{r\in I}}\binom{j}{%
s_{1}s_{2...}s_{n}}\Pi _{r=1}^{n}c_{r}^{s_{r}}\Psi
(k-\sum_{r=1}^{n}(m+1)s_{r}+\sum_{r=1}^{n}rs_{r})x^{k}y^{j} \\
&=&\sum_{j=0}^{\infty }\sum_{k=-\infty }^{\infty }\sum_{s_{r\in I}}\binom{j}{%
s_{1}s_{2...}s_{n}}\Pi _{r=1}^{n}c_{r}^{s_{r}}\Psi
(k+\sum_{r=1}^{n}rs_{r}-(m+1)j)x^{k}y^{j}.
\end{eqnarray*}
Thus the solution is of the form 
\[
U_{i,j}=\sum_{s_{r\in I}}\binom{j}{s_{1}s_{2...}s_{n}}\Pi
_{r=1}^{n}c_{r}^{s_{r}}\Psi (i+\sum_{r=1}^{n}rs_{r}-(m+1)j), 
\]
where $I=\{s_{1},s_{2},...,s_{n}|0\leq s_{k}\leq j$, $\sum_{k=1}^{n}s_{k}=j%
\}.$%

\endproof%
Next, we turn to (2.3). Suppose $U_{i,j}$ is known on x-axis, in general (2.3) has infinitely many solutions. Consider the typical example where $a=1$, $b=1$, and $c=-1$ then $U_{i,j}=c_{i,j}$ is a solution. Since there are infinitely many choices of constants $c_{j}$ for $j=0,1,...$, we then have infinite solutions. Using generating functions, we can find one special solution of (2.3).

\begin{theorem}
Let 
\[
U_{i+1,j+1}=aU_{i,j+1}+bU_{i+1,j}+cU_{i,j} 
\]
with i$\in Z$, j$\in N\cup \{0\}$
Given $U_{i,0}=\Psi (i)$, then it has following solution 
\[
U_{i,j}=\sum_{s=0}^{\infty }\Omega (i-s)p(s,j), 
\]
where 
\begin{eqnarray*}
\Omega (k) &=&\Psi (k)-a\Psi (k-1), \\
p(j,t) &=&\sum_{r=0}^{t}\binom{j}{t}\binom{t}{r}a^{j-t}b^{t-r}c^{t}.
\end{eqnarray*}

\proof%
Let 
\[
\sum_{i=-\infty }^{\infty }\sum_{j=0}^{\infty }U_{i,j}x^{i}y^{j}=h(x,y). 
\]
Then we have
\end{theorem}

\begin{eqnarray*}
&&\sum_{i=-\infty }^{\infty }\sum_{j=0}^{\infty }U_{i+1,j+1}x^{i}y^{j} \\
&=&a\sum_{i=-\infty }^{\infty }\sum_{j=0}^{\infty
}U_{i,j+1}x^{i}y^{j}+b\sum_{i=-\infty }^{\infty }\sum_{j=0}^{\infty
}U_{i+1,j}x^{i}y^{j}+c\sum_{i=-\infty }^{\infty }\sum_{j=0}^{\infty
}U_{i,j}x^{i}y^{j},
\end{eqnarray*}

\begin{eqnarray*}
\frac{1}{xy}[h(x,y)-\sum_{i=-\infty }^{\infty }\Psi (i)x^{i}] &=&\frac{a}{y}[%
h(x,y)-\sum_{i=-\infty }^{\infty }\Psi (i)x^{i}]+\frac{b}{x}h(x,y)+ch(x,y),
\\
h(x,y)[1-(ax+by+cxy)] &=&\sum_{i=-\infty }^{\infty }\Psi
(i)x^{i}-a\sum_{i=-\infty }^{\infty }\Psi (i)x^{i+1} \\
&=&\sum_{k=-\infty }^{\infty }[\Psi (k)-a\Psi (k-1)]x^{k},
\end{eqnarray*}

\begin{eqnarray*}
h(x,y) &=&\frac{\sum_{k=-\infty }^{\infty }[\Psi (k)-a\Psi (k-1)]x^{k}}{%
1-(ax+by+cxy)}{} \\
&=&\sum_{k=-\infty }^{\infty }[\Psi (k)-a\Psi (k-1)]x^{k}\sum_{j=0}^{\infty
}(ax+by+cxy)^{j} \\
&=&\sum_{k=-\infty }^{\infty }\sum_{j=0}^{\infty }[\Psi (k)-a\Psi
(k-1)]x^{k}(ax+(b+cx)y)^{j} \\
&=&\sum_{k=-\infty }^{\infty }\sum_{j=0}^{\infty }[\Psi (k)-a\Psi
(k-1)]x^{k}\sum_{t=0}^{j}\binom{j}{t}(ax)^{j-t}(b+cx)^{t}y^{t} \\
&=&\sum_{k=-\infty }^{\infty }\sum_{j=0}^{\infty }[\Psi (k)-a\Psi
(k-1)]x^{k}\sum_{t=0}^{j}\binom{j}{t}[\sum_{r=0}^{t}\binom{t}{r}%
b^{t-r}c^{t}x^{t}]a^{j-t}x^{j-t}y^{t} \\
&=&\sum_{k=-\infty }^{\infty }\sum_{j=0}^{\infty }[\Psi (k)-a\Psi
(k-1)]x^{k}\sum_{t=0}^{j}\sum_{r=0}^{t}\binom{j}{t}\binom{t}{r}%
b^{t-r}c^{t}x^{t}a^{j-t}x^{j-t}y^{t} \\
&=&\sum_{k=-\infty }^{\infty }\sum_{j=0}^{\infty }[\Psi (k)-a\Psi
(k-1)]x^{k}\sum_{t=0}^{j}\sum_{r=0}^{t}\binom{j}{t}\binom{t}{r}%
a^{j-t}b^{t-r}c^{t}x^{j}y^{t} \\
&=&\sum_{k=-\infty }^{\infty }\sum_{j=0}^{\infty }[\Psi (k)-a\Psi
(k-1)]\sum_{t=0}^{j}\sum_{r=0}^{t}\binom{j}{t}\binom{t}{r}%
a^{j-t}b^{t-r}c^{t}x^{j+k}y^{t}.
\end{eqnarray*}
Let 
\begin{eqnarray*}
\Psi (k)-a\Psi (k-1) &=&\Omega (k), \\
\sum_{r=0}^{t}\binom{j}{t}\binom{t}{r}a^{j-t}b^{t-r}c^{t} &=&p(j,t), \\
j+k &=&u.
\end{eqnarray*}
Then we have 
\begin{eqnarray*}
h(x,y) &=&\sum_{u=-\infty }^{\infty }\sum_{j=0}^{\infty }\Omega
(u-j)\sum_{t=0}^{j}p(j,t)x^{u}y^{t} \\
&=&\sum_{u=-\infty }^{\infty }x^{u}\sum_{j=0}^{\infty }\Omega
(u-j)\sum_{t=0}^{j}p(j,t)y^{t} \\
&=&\sum_{u=-\infty }^{\infty }x^{u}\Omega (u)p(0,0)+ \\
&&\sum_{u=-\infty }^{\infty }x^{u}\Omega (u-1)[p(1,0)+p(1,1)y]+...+ \\
&&\sum_{u=-\infty }^{\infty }x^{u}\Omega
(u-k)[p(k,0)+p(k,1)y+...+p(k,k)y^{k}]+... \\
&=&\sum_{u=-\infty }^{\infty }x^{u}\sum_{j=0}^{\infty }\sum_{s=0}^{\infty
}\Omega (u-s)p(s,j)y^{j} \\
&=&\sum_{i=-\infty }^{\infty }\sum_{j=0}^{\infty }\sum_{s=0}^{\infty }\Omega
(i-s)p(s,j)x^{i}y^{j}
\end{eqnarray*}

Thus the solution is 
\[
U_{i,j}=\sum_{s=0}^{\infty }\Omega (i-s)p(s,j), 
\]
where 
\begin{eqnarray*}
\Omega (k) &=&\Psi (k)-a\Psi (k-1), \\
p(j,t) &=&\sum_{r=0}^{t}\binom{j}{t}\binom{t}{r}a^{j-t}b^{t-r}c^{t}.
\end{eqnarray*}
\endproof%
Notice that in above theorem, we have the initial condition on x-axis and then obtain the solution on the upper half plane. Suppose we are only given the initial data on an interval or ray of x-axis. On its domain of influence, the solutions of (\ref{eq1}), (\ref{eq2}), (\ref{eq3}) are also given by the corresponding theorems.

\section{The General Case}
In this section, we discuss the general solution of the partial difference
equation 
\begin{equation}
U_{i+m,j+k}=\sum_{s=1}^{k}\sum_{r=1}^{n}c_{rs}U_{i+(r-1),j+(s-1)},
\label{eq5}
\end{equation}
which is an extension of the equation(\ref{eq2}) in the previous section. We try to derive its explicit solution with given initial conditions data $U_{i,t}=\Psi
_{t}(i)$, $t=1,2,...,k-1$.The theorem below is the case for $k=2$.

\begin{theorem}
Let 
\[
U_{i+m,j+k}=\sum_{s=1}^{k}\sum_{r=1}^{n}c_{rs}U_{i+(r-1),j+(s-1)} 
\]
\end{theorem}
with i$\in Z$, $j\in N\cup\{0\} $. 
Given $U_{i,t}=\Psi _{t}(i)$,for $t=0,1$, then its solution $U_{i,j}$ is of the form

\begin{eqnarray*}
&&\sum_{l_{rs}\in I_{1}}\binom{\frac{j+\sum_{s=1}^{2}\sum_{r=1}^{n}sl_{rs}}{3%
}}{l_{11}l_{12}...l_{n2}}\prod_{\substack{s=1,2,\\r=1,2,...n}}c_{rs}^{lrs}
\\
&&\Psi _{0}(i+\sum_{s=1}^{2}\sum_{r=1}^{n}rl_{rs}-(m+1)\frac{%
j+\sum_{s=1}^{2}\sum_{r=1}^{n}sl_{rs}}{3})+ \\
&&\sum_{l_{rs}\in I_{2}}\binom{\frac{(j-1+\sum_{s=1}^{2}%
\sum_{r=1}^{n}sl_{rs})}{3}}{l_{11}l_{12}...l_{n2}}\prod_{\substack{s=1,2,\\r=1,2,...n}}c_{rs}^{lrs}\Psi
_{1}(i+\sum_{s=1}^{2}\sum_{r=1}^{n}rl_{rs}-(m+1) \\
&&\frac{(j-1+\sum_{s=1}^{2}\sum_{r=1}^{n}sl_{rs})}{3})- \\
&&\sum_{e=1}^{n}c_{e2}\sum_{l_{rs}\in I_{2}}\binom{\frac{(j-1+\sum_{s=1}^{2}%
\sum_{r=1}^{n}sl_{rs})}{3}}{l_{11}l_{12}...l_{n2}}\prod_{\substack{s=1,2,\\r=1,2,...n}}c_{rs}^{lrs}\\
&&\Psi _{0}(i+e-1-m+\sum_{s=1}^{2}\sum_{r=1}^{n}rl_{rs}-(m+1)\frac{%
(j-1+\sum_{s=1}^{2}\sum_{r=1}^{n}sl_{rs})}{3}),
\end{eqnarray*}
where 
\begin{eqnarray*}
I_{1} &=&\{l_{11},l_{12},...,l_{n2}|\text{ }0\leq l_{rs}\leq \frac{%
j+\sum_{s=1}^{2}\sum_{r=1}^{n}sl_{rs}}{3}, \\
\sum_{s=1}^{2}\sum_{r=1}^{n}l_{rs} &=&\frac{j+\sum_{s=1}^{2}%
\sum_{r=1}^{n}sl_{rs}}{3}\}, \\
I_{2} &=&\{l_{11},l_{12},...,l_{n2}|\text{ }0\leq l_{rs}\leq \frac{%
(j-1+\sum_{s=1}^{2}\sum_{r=1}^{n}sl_{rs})}{3}, \\
\sum_{s=1}^{2}\sum_{r=1}^{n}l_{rs} &=&\frac{(j-1+\sum_{s=1}^{2}%
\sum_{r=1}^{n}sl_{rs})}{3}\},
\end{eqnarray*}
for s$=$1, 2, r=1, 2, ..., n.

\proof%
Let $h(x,y)=\sum_{i=-\infty }^{\infty }\sum_{j=0}^{\infty }U_{i,j}x^{i}y^{j}$%
and $k=2$. Then we have

\[
\sum_{i=-\infty }^{\infty }\sum_{j=0}^{\infty
}U_{i+m,j+k}x^{i}y^{j}=\sum_{s=1}^{k}%
\sum_{r=1}^{n}c_{rs}U_{i+(r-1),j+(s-1)}x^{i}y^{j}, 
\]
\begin{eqnarray*}
&&\frac{1}{x^{m}y^{k}}[h(x,y)-\sum_{r=0}^{k-1}\sum_{i=-\infty }^{\infty
}\Psi _{r}(i)x^{i}y^{r}] \\
&=&\sum_{r=1}^{n}c_{r1}x^{-(r-1)}h(x,y)+%
\sum_{r=1}^{n}c_{r2}x^{-(r-1)}y^{-1}[h(x,y)-\sum_{i=-\infty }^{\infty }\Psi
_{0}(i)x^{i}]+ \\
&&\sum_{r=1}^{n}c_{r3}x^{-(r-1)}y^{-2}[h(x,y)-\sum_{t=0}^{1}\sum_{i=-\infty
}^{\infty }\Psi _{t}(i)x^{i}y^{t}]+...+ \\
&&\sum_{r=1}^{n}c_{rk}x^{-(r-1)}y^{-(k-1)}[h(x,y)-\sum_{t=0}^{k-2}\sum_{i=-%
\infty }^{\infty }\Psi _{t}(i)x^{i}y^{t}],
\end{eqnarray*}
\begin{eqnarray*}
&&L(x,y)\\
&=&h(x,y)[1-(\sum_{r=1}^{n}c_{r1}x^{m-(r-1)}y^{k}+%
\sum_{r=1}^{n}c_{r2}x^{m-(r-1)}y^{k-1}+...+\sum_{r=1}^{n}c_{rk}x^{m-(r-1)}y) \\
&=&\sum_{r=0}^{k-1}\sum_{i=-\infty }^{\infty }\Psi
_{r}(i)x^{i}y^{r}-\sum_{s=2}^{k}[\sum_{r=1}^{n}c_{rs}x^{m-(r-1)}y^{k-s+1}%
\sum_{t=0}^{s-2}\sum_{i=-\infty }^{\infty }\Psi _{t}(i)x^{i}y^{t}] \\
&=&\sum_{i=-\infty }^{\infty }\sum_{r=0}^{k-1}\Psi
_{r}(i)x^{i}y^{r}-\sum_{i=-\infty }^{\infty
}\sum_{s=2}^{k}\sum_{t=0}^{s-2}\sum_{r=1}^{n}c_{rs}\Psi
_{t}(i)x^{m-(r-1)+i}y^{k-s+1+t} \\
&=&\sum_{i=-\infty }^{\infty }\sum_{r=0}^{k-1}\Psi
_{r}(i)x^{i}y^{r}-\sum_{i=-\infty }^{\infty
}\sum_{h=2}^{k}\sum_{t=0}^{h-2}\sum_{e=1}^{n}c_{eh}\Psi
_{t}(i)x^{i+m-e+1}y^{k-h+1+t}, \\
\end{eqnarray*}
\begin{eqnarray*}
&&h(x,y)\\
&=&L(x,y)[1-\sum_{s=1}^{k}\sum_{r=1}^{n}c_{rs}x^{m-(r-1)}y^{k-s+1}]^{-1} \\
&=&L(x,y)\sum_{j=0}^{\infty
}[\sum_{s=1}^{k}\sum_{r=1}^{n}c_{rs}x^{m-(r-1)}y^{k-s+1}]^{j} \\
&=&L(x,y)\sum_{j=0}^{\infty }\sum_{l_{rs}\in I}\binom{j}{%
l_{11}l_{12}...l_{1n}l_{21}...l_{nk}}\prod_{\substack{r=1,2,...,n,\\s=1,2,...,k}} 
c_{rs}^{l_{rs}}x^{a}y^{b} \\
&=&\sum_{j=0}^{\infty }L(x,y)\sum_{l_{rs}\in I}\binom{j}{%
l_{11}l_{12}...l_{1n}l_{21}...l_{nk}}\prod_{\substack{r=1,2,...,n,\\s=1,2}} 
c_{rs}^{l_{rs}}x^{a}y^{b},
\end{eqnarray*}
where 
\begin{eqnarray*}
I &=&\{l_{11},l_{12},...,l_{nk}|\text{ }0\leq l_{rs}\leq j\text{, }%
\sum_{s=1}^{k}\sum_{r=1}^{n}l_{rs}=j\}, \\
a &=&\sum_{s=1}^{k}\sum_{r=1}^{n}(m-r+1)l_{rs} \\
&=&(m+1)j-\sum_{s=1}^{k}\sum_{r=1}^{n}rl_{rs}\text{,} \\
b &=&\sum_{s=1}^{k}\sum_{r=1}^{n}(k-s+1)l_{rs} \\
&=&(k+1)\sum_{s=1}^{k}\sum_{r=1}^{n}l_{rs}-\sum_{s=1}^{k}%
\sum_{r=1}^{n}sl_{rs} \\
&=&(k+1)j-\sum_{s=1}^{k}\sum_{r=1}^{n}sl_{rs},
\end{eqnarray*}
for r=1, 2, ..., n, s=1, 2, ..., k. When k=2,we have 
\begin{eqnarray*}
I &=&\{l_{11},l_{12},...,l_{n\times 2}|\text{ }0\leq l_{rs}\leq j\text{, }%
\sum_{s=1}^{2}\sum_{r=1}^{n}l_{rs}=j\}, \\
a &=&(m+1)j-\sum_{s=1}^{2}\sum_{r=1}^{n}rl_{rs}\text{,} \\
b &=&3j-\sum_{s=1}^{2}\sum_{r=1}^{n}sl_{rs},
\end{eqnarray*}
and 
\begin{eqnarray*}
&&L(x,y)=\sum_{i=-\infty }^{\infty }\sum_{r=0}^{1}\Psi
_{r}(i)x^{i}y^{r}-\sum_{i=-\infty }^{\infty }\sum_{e=1}^{n}c_{e2}\Psi
_{0}(i)x^{i+m-e+1}y \\
&=&\sum_{i=-\infty }^{\infty }\Psi _{0}(i)x^{i}+\sum_{i=-\infty }^{\infty
}\Psi _{1}(i)x^{i}y-\sum_{i=-\infty }^{\infty }\sum_{e=1}^{n}c_{e2}\Psi
_{0}(i)x^{i+m-e+1}y
\end{eqnarray*}
Hence
\begin{eqnarray*}
&&h(x,y)=\sum_{j=0}^{\infty }\{\sum_{i=-\infty }^{\infty }\Psi
_{0}(i)x^{i}+\sum_{i=-\infty }^{\infty }\Psi _{1}(i)x^{i}y-\sum_{i=-\infty
}^{\infty }\sum_{e=1}^{n}c_{e2}\Psi _{0}(i)x^{i+m-e+1}y\} \\
&&\sum_{l_{rs}\in I}\binom{j}{l_{11}l_{12}...l_{1n}l_{21}...l_{nk}}\prod_{\substack{r=1,2,...,n,\\s=1,2,...,k}} c_{rs}^{l_{rs}}x^{a}y^{b} \\
&=&\sum_{j=0}^{\infty }\sum_{i=-\infty }^{\infty }\Psi
_{0}(i)x^{i}\sum_{l_{rs}\in I}\binom{j}{l_{11}l_{12}...l_{1n}l_{21}...l_{nk}}%
\prod_{\substack{r=1,2,...,n,\\s=1,2,...,k}}c_{rs}^{l_{rs}}x^{a}y^{b}+ \\
&&\sum_{j=0}^{\infty }\sum_{i=-\infty }^{\infty }\Psi
_{1}(i)x^{i}y\sum_{l_{rs}\in I}\binom{j}{l_{11}l_{12}...l_{1n}l_{21}...l_{nk}%
}\prod_{\substack{r=1,2,...,n,\\s=1,2,...,k}}c_{rs}^{l_{rs}}x^{a}y^{b}- \\
&&\sum_{j=0}^{\infty }\sum_{i=-\infty }^{\infty }[\sum_{e=1}^{n}c_{e2}\Psi
_{0}(i)x^{i+m-e+1}]y\sum_{l_{rs}\in I}\binom{j}{%
l_{11}l_{12}...l_{1n}l_{21}...l_{nk}}\prod_{\substack{r=1,2,...,n,\\s=1,2,...,k}}c_{rs}^{l_{rs}}x^{a}y^{b} \\
&=&\sum_{j=0}^{\infty }\sum_{i=-\infty }^{\infty }\sum_{l_{rs}\in I}\binom{j%
}{l_{11}l_{12}...l_{1n}l_{21}...l_{nk}}\Psi _{0}(i)\prod_{r=1,2,...,n, s=1,2,...,k}c_{rs}^{l_{rs}}x^{a+i}y^{b}+ \\
&&\sum_{j=0}^{\infty }\sum_{i=-\infty }^{\infty }\sum_{l_{rs}\in I}\binom{j}{%
l_{11}l_{12}...l_{1n}l_{21}...l_{nk}}\Psi _{1}(i)\prod_{\substack{r=1,2,...,n,\\s=1,2,...,k}}c_{rs}^{l_{rs}}x^{a+i}y^{b+1}- \\
&&\sum_{j=0}^{\infty }\sum_{i=-\infty }^{\infty
}\sum_{e=1}^{n}c_{e2}\sum_{l_{rs}\in I}\binom{j}{%
l_{11}l_{12}...l_{1n}l_{21}...l_{nk}}\Psi _{0}(i)\prod_{\substack{r=1,2,...,n,\\s=1,2,...,k}} c_{rs}^{l_{rs}}x^{a}y^{b}x^{a+i+m-e+1}y^{b+1}
\end{eqnarray*}
Let 
\[
a+i=\alpha ,b+1=3j-\sum_{s=1}^{2}\sum_{r=1}^{n}sl_{rs}+1=\beta
,a+i+m-e+1=\gamma . 
\]
Then

\begin{eqnarray*}
j &=&\frac{b+\sum_{s=1}^{2}\sum_{r=1}^{n}sl_{rs}}{3}=\frac{(\beta
-1+\sum_{s=1}^{2}\sum_{r=1}^{n}sl_{rs})}{3}, \\
i &=&\alpha -a \\
&=&\alpha -(m+1)j+\sum_{s=1}^{2}\sum_{r=1}^{n}rl_{rs} \\
&=&\alpha +\sum_{s=1}^{2}\sum_{r=1}^{n}rl_{rs}-(m+1)\frac{%
b+\sum_{s=1}^{2}\sum_{r=1}^{n}sl_{rs}}{3} \\
&=&\alpha +\sum_{s=1}^{2}\sum_{r=1}^{n}rl_{rs}-(m+1)\frac{(\beta
-1+\sum_{s=1}^{2}\sum_{r=1}^{n}sl_{rs})}{3} \\
&=&\gamma +e-1-m-a \\
&=&\gamma +e-1-m+\sum_{s=1}^{2}\sum_{r=1}^{n}rl_{rs}-(m+1)\frac{(\beta
-1+\sum_{s=1}^{2}\sum_{r=1}^{n}sl_{rs})}{3},
\end{eqnarray*}
and
\begin{eqnarray*}
&&h(x,y)=\sum_{b=0}^{\infty }\sum_{\alpha =-\infty }^{\infty
}\sum_{l_{rs}\in I_{1}}\binom{\frac{b+\sum_{s=1}^{2}\sum_{r=1}^{n}sl_{rs}}{3}%
}{l_{11}l_{12}...l_{n2}}\prod_{\substack{s=1,2,\\r=1,2,...n}}%
c_{rs}^{lrs}x^{\alpha }y^{b} \\
&&\Psi _{0}(\alpha +\sum_{s=1}^{2}\sum_{r=1}^{n}rl_{rs}-(m+1)\frac{%
b+\sum_{s=1}^{2}\sum_{r=1}^{n}sl_{rs}}{3})+ \\
&&\sum_{\beta =0}^{\infty }\sum_{\alpha =-\infty }^{\infty }\sum_{l_{rs}\in
I_{2}}\binom{\frac{(\beta -1+\sum_{s=1}^{2}\sum_{r=1}^{n}sl_{rs})}{3}}{%
l_{11}l_{12}...l_{n2}}\prod_{\substack{s=1,2,\\r=1,2,...n}}%
c_{rs}^{lrs}x^{\alpha }y^{\beta } \\
&&\Psi _{1}(\alpha +\sum_{s=1}^{2}\sum_{r=1}^{n}rl_{rs}-(m+1)\frac{(\beta
-1+\sum_{s=1}^{2}\sum_{r=1}^{n}sl_{rs})}{3})- \\
&&\sum_{\beta =0}^{\infty }\sum_{\gamma =-\infty }^{\infty
}\sum_{e=1}^{n}c_{e2}\sum_{l_{rs}\in I_{2}}\binom{\frac{(\beta
-1+\sum_{s=1}^{2}\sum_{r=1}^{n}sl_{rs})}{3}}{l_{11}l_{12}...l_{n2}}\prod_{\substack{s=1,2,\\r=1,2,...n}}%
c_{rs}^{lrs}x^{\gamma }y^{\beta } \\
&&\Psi _{0}(\gamma +e-1-m+\sum_{s=1}^{2}\sum_{r=1}^{n}rl_{rs}-(m+1)\frac{%
(\beta -1+\sum_{s=1}^{2}\sum_{r=1}^{n}sl_{rs})}{3}) \\
&=&\sum_{u=0}^{\infty }\sum_{v=-\infty }^{\infty }\sum_{l_{rs}\in I_{1}}%
\binom{\frac{u+\sum_{s=1}^{2}\sum_{r=1}^{n}sl_{rs}}{3}}{l_{11}l_{12}...l_{n2}%
}\prod_{\substack{s=1,2,\\r=1,2,...n}}c_{rs}^{lrs} \\
&&\Psi _{0}(v+\sum_{s=1}^{2}\sum_{r=1}^{n}rl_{rs}-(m+1)\frac{%
u+\sum_{s=1}^{2}\sum_{r=1}^{n}sl_{rs}}{3})+ \\
&&\sum_{l_{rs}\in I_{2}}\binom{\frac{(u-1+\sum_{s=1}^{2}%
\sum_{r=1}^{n}sl_{rs})}{3}}{l_{11}l_{12}...l_{n2}}\prod_{\substack{s=1,2,\\r=1,2,...n}}c_{rs}^{lrs}\Psi
_{1}(v+\sum_{s=1}^{2}\sum_{r=1}^{n}rl_{rs}-(m+1) \\
&&\frac{(u-1+\sum_{s=1}^{2}\sum_{r=1}^{n}sl_{rs})}{3})- \\
&&\sum_{e=1}^{n}c_{e2}\sum_{l_{rs}\in I_{2}}\binom{\frac{(u-1+\sum_{s=1}^{2}%
\sum_{r=1}^{n}sl_{rs})}{3}}{l_{11}l_{12}...l_{n2}}\prod_{\substack{s=1,2,\\r=1,2,...n}}c_{rs}^{lrs} \\
&&\Psi _{0}(v+e-1-m+\sum_{s=1}^{2}\sum_{r=1}^{n}rl_{rs}-(m+1)\frac{%
(u-1+\sum_{s=1}^{2}\sum_{r=1}^{n}sl_{rs})}{3})x^{v}y^{u}.
\end{eqnarray*}

Thus for k=2, the solution $U_{i,j}$ has the desired form.

\endproof%
\section{The Multidimensional Cases}
In this section, we discuss the linear partial difference equations in higher dimensions. At the beging, we solve the following equations in three dimensions where $i$, $j\in Z$, and $k\in N\cup \{0\} $.

\begin{equation} 
\begin{split}
U_{i,j,k+1}=c_{1}U_{i-1,j-1,k}+c_{2}U_{i,j-1,k}+c_{3}U_{i+1,j-1,k}+c_{4}U_{i-1,j,k}+\\
c_{5}U_{i,j,k}+c_{6}U_{i+1,j,k}+c_{7}U_{i-1,j+1,k}+c_{8}U_{i,j+1,k}+c_{9}U_{i+1,j+1,k}
\label{eq6}
\end{split}
\end{equation}

\begin{equation}
U_{i+s,j+t,k+1}=\sum_{u=1}^{n}\sum_{v=1}^{m}c_{uv}U_{i+u-1,,j+v-1,k}
\label{eq7}
\end{equation}
Note that equation(\ref{eq1}) in section 2 is a special case for equation(\ref{eq6}) and equation(\ref{eq2}) is a special case for equation(\ref{eq7}). We will obtain the explicit solution under the assumption of given $U_{i,j,0}=\Psi (i,j).$

\begin{theorem}
Let 
\begin{eqnarray*}
U_{i,j,k+1} &=&c_{1}U_{i-1,j-1,k}+c_{2}U_{i,j-1,k}+c_{3}U_{i+1,j-1,k}+ \\
&&c_{4}U_{i-1,j,k}+c_{5}U_{i,j,k}+c_{6}U_{i+1,j,k}+ \\
&&c_{7}U_{i-1,j+1,k}+c_{8}U_{i,j+1,k}+c_{9}U_{i+1,j+1,k}
\end{eqnarray*}

with i, j$\in Z$, $%
k\in N\cup \{0\}$. \ Given $U_{i,j,0}=\Psi (i,j)$, then its solution is of
the form

\[
U_{i,j,k}=\sum_{s_{r\in I}}\binom{k}{s_{1},s_{2},...,s_{9}}\Psi
(i-a,j-b)\prod_{r=1}^{9}c_{r}^{sr}, 
\]
where $I=\{s_{1},s_{2},...,s_{9}|0\leq s_{r}\leq k$, $\sum_{r=1}^{9}s_{r}=k%
\} $, $a=s_{1}-s_{3}+s_{4}-s_{6}+s_{7}-s_{9}$, and $%
b=s_{1}+s_{2}+s_{3}-s_{7}-s_{8}-s_{9}$.
\end{theorem}
\proof%
Let 
\[
h(x,y,z)=\sum_{i=-\infty }^{\infty }\sum_{j=-\infty }^{\infty
}\sum_{k=0}^{\infty }U_{i,j,k}x^{i}y^{j}z^{k}. 
\]
Then we have

\begin{eqnarray*}
\sum_{i=-\infty }^{\infty }\sum_{j=-\infty }^{\infty }\sum_{k=0}^{\infty
}U_{i,j,k+1}x^{i}y^{j}z^{k} &=&c_{1}\sum_{i=-\infty }^{\infty
}\sum_{j=-\infty }^{\infty }\sum_{k=0}^{\infty }U_{i-1,j-1,k}x^{i}y^{j}z^{k}+
\\
&&c_{2}\sum_{i=-\infty }^{\infty }\sum_{j=-\infty }^{\infty
}\sum_{k=0}^{\infty }U_{i,j-1,k}x^{i}y^{j}z^{k}+ \\
&&c_{3}\sum_{i=-\infty }^{\infty }\sum_{j=-\infty }^{\infty
}\sum_{k=0}^{\infty }U_{i+1,j-1,k}x^{i}y^{j}z^{k}+ \\
&&c_{4}\sum_{i=-\infty }^{\infty }\sum_{j=-\infty }^{\infty
}\sum_{k=0}^{\infty }U_{i-1,j,k}x^{i}y^{j}z^{k}+ \\
&&c_{5}\sum_{i=-\infty }^{\infty }\sum_{j=-\infty }^{\infty
}\sum_{k=0}^{\infty }U_{i,j,k}x^{i}y^{j}z^{k}+ \\
&&c_{6}\sum_{i=-\infty }^{\infty }\sum_{j=-\infty }^{\infty
}\sum_{k=0}^{\infty }U_{i+1,j,k}x^{i}y^{j}z^{k}+ \\
&&c_{7}\sum_{i=-\infty }^{\infty }\sum_{j=-\infty }^{\infty
}\sum_{k=0}^{\infty }U_{i-1,j+1,k}x^{i}y^{j}z^{k}+ \\
&&c_{8}\sum_{i=-\infty }^{\infty }\sum_{j=-\infty }^{\infty
}\sum_{k=0}^{\infty }U_{i,j+1,k}x^{i}y^{j}z^{k}+ \\
&&c_{9}\sum_{i=-\infty }^{\infty }\sum_{j=-\infty }^{\infty
}\sum_{k=0}^{\infty }U_{i+1,j+1,k}x^{i}y^{j}z^{k},
\end{eqnarray*}

\begin{eqnarray*}
\frac{1}{z}[h(x,y,z)-\sum_{i=-\infty }^{\infty }\sum_{j=-\infty }^{\infty
}\Psi (i,j)x^{i}y^{j}]
&=&(c_{1}xy+c_{2}y+c_{3}x^{-1}y+c_{4}x+c_{5}+c_{6}x^{-1}+ \\
&&c_{7}xy^{-1}+c_{8}y^{-1}+c_{9}x^{-1}y^{-1})h(x,y,z),
\end{eqnarray*}

\begin{eqnarray*}
\sum_{i=-\infty }^{\infty }\sum_{j=-\infty }^{\infty }\Psi (i,j)x^{i}y^{j}
&=&[1-(c_{1}xy+c_{2}y+c_{3}x^{-1}y+c_{4}x+c_{5}+ \\
&&c_{6}x^{-1}+c_{7}xy^{-1}+c_{8}y^{-1}+c_{9}x^{-1}y^{-1})z]h(x,y,z),
\end{eqnarray*}
\begin{eqnarray*}
h(x,y,z) &=&\frac{\sum_{i=-\infty }^{\infty }\sum_{j=-\infty }^{\infty }\Psi
(i,j)x^{i}y^{j}}{%
1-(c_{1}xy+c_{2}y+c_{3}x^{-1}y+c_{4}x+c_{5}+c_{6}x^{-1}+c_{7}xy^{-1}+c_{8}y^{-1}+c_{9}x^{-1}y^{-1})z%
} \\
&=&\sum_{i=-\infty }^{\infty }\sum_{j=-\infty }^{\infty }\Psi
(i,j)x^{i}y^{j}\sum_{k=0}^{\infty }(c_{1}xy+c_{2}y+c_{3}x^{-1}y+c_{4}x+c_{5}+
\\
&&c_{6}x^{-1}+c_{7}xy^{-1}+c_{8}y^{-1}+c_{9}x^{-1}y^{-1})^{k}z^{k} \\
&=&\sum_{k=0}^{\infty }\sum_{i=-\infty }^{\infty }\sum_{j=-\infty }^{\infty
}\Psi (i,j)x^{i}y^{j}\sum_{s_{r\in I}}\binom{k}{s_{1},s_{2},...,s_{9}}%
c_{1}^{s_{1}}(xy)^{s_{1}}c_{2}^{s_{2}}y^{s_{2}}...c_{9}^{s_{9}}(xy)^{-s_{9}}z^{k}
\\
&=&\sum_{k=0}^{\infty }\sum_{i=-\infty }^{\infty }\sum_{j=-\infty }^{\infty
}\Psi (i,j)x^{i}y^{j}\sum_{s_{r\in I}}\binom{k}{s_{1},s_{2},...,s_{9}}%
\prod_{r=1}^{9}c_{r}^{sr}x^{a}y^{b}z^{k} \\
&=&\sum_{k=0}^{\infty }\sum_{i=-\infty }^{\infty }\sum_{j=-\infty }^{\infty
}\sum_{s_{r\in I}}\binom{k}{s_{1},s_{2},...,s_{9}}\Psi
(i,j)\prod_{r=1}^{9}c_{r}^{sr}x^{a+i}y^{b+j}z^{k},
\end{eqnarray*}
where $I=\{s_{1},s_{2},...,s_{9}|0\leq s_{r}\leq k$, $\sum_{r=1}^{9}s_{r}=k%
\} $, $a=s_{1}-s_{3}+s_{4}-s_{6}+s_{7}-s_{9}$, and $%
b=s_{1}+s_{2}+s_{3}-s_{7}-s_{8}-s_{9}$. Let 
\[
a+i=u\text{, }b+j=v. 
\]
Then we have 
\[
h(x,y,z)=\sum_{k=0}^{\infty }\sum_{u=-\infty }^{\infty }\sum_{v=-\infty
}^{\infty }\sum_{s_{r\in I}}\binom{k}{s_{1},s_{2},...,s_{9}}\Psi
(u-a,v-b)\prod_{r=1}^{9}c_{r}^{sr}x^{u}y^{v}z^{k}. 
\]
Hence the solution is as follows 
\[
U_{i,j,k}=\sum_{s_{r\in I}}\binom{k}{s_{1},s_{2},...,s_{9}}\Psi
(i-a,j-b)\prod_{r=1}^{9}c_{r}^{sr}, 
\]
\endproof%
\begin{theorem}
Let 
\[
U_{i+s,j+t,k+1}=\sum_{u=1}^{n}\sum_{v=1}^{m}c_{uv}U_{i+u-1,,j+v-1,k} 
\]
with i$\in Z$, $j\in Z$, $k\in
N\cup \{0\}$. \ Given $U_{i,j,0}=\Psi (i,j)$, then its solution is of the form

\[
U_{i,j}=\sum_{s_{qr\in I}}\binom{k}{s_{11},s_{12},...,s_{1m},...,s_{nm}}\Psi
(i-a,j-b)\prod_{\substack{u=1,2,...,n,\\v=1,2,...,m}}c_{uv}^{s_{uv}}, 
\]
where $I=\{s_{11},s_{12},...,s_{1m},...,s_{nm}|0\leq s_{uv}\leq
j-b,\sum_{u=1}^{n}\sum_{v=1}^{m}s_{uv}=j-b\}$, $a=(s+1)k-\sum_{u=1}^{n}%
\sum_{v=1}^{m}us_{uv}$, and $b=(t+1)k-\sum_{u=1}^{n}\sum_{v=1}^{m}vs_{uv}.$
\end{theorem}
Notice that this theorem could provide a mathematical framework for understanding how information within a black hole is encoded on its event horizon, a key concept in the black hole information paradox.  
\proof%
Let 
\[
h(x,y,z)=\sum_{i=-\infty }^{\infty }\sum_{j=-\infty }^{\infty
}\sum_{k=0}^{\infty }U_{i,j}x^{i}y^{j}z^{k}. 
\]
Then we have

\[
\sum_{i=-\infty }^{\infty }\sum_{j=-\infty }^{\infty }\sum_{k=0}^{\infty
}U_{i+s,j+t,k+1}x^{i}y^{j}z^{k}=\sum_{i=-\infty }^{\infty }\sum_{j=-\infty
}^{\infty }\sum_{k=0}^{\infty
}\sum_{u=1}^{n}\sum_{v=1}^{m}c_{uv}U_{i+u-1,,j+v-1,k}x^{i}y^{j}z^{k}, 
\]

\[
\frac{1}{x^{s}y^{t}z}[h(x,y,z)-\sum_{i=-\infty }^{\infty }\sum_{j=-\infty
}^{\infty }\Psi (i,j)x^{i}y^{j}]=\sum_{u=1}^{n}%
\sum_{v=1}^{m}c_{uv}x^{-(u-1)}y^{-(v-1)}h(x,y,z), 
\]

\[
h(x,y)(1-\sum_{u=1}^{n}\sum_{v=1}^{m}c_{uv}x^{s-(u-1)}y^{t-(v-1)}z)=%
\sum_{i=-\infty }^{\infty }\sum_{j=-\infty }^{\infty }\Psi (i,j)x^{i}y^{j}, 
\]

\begin{eqnarray*}
h(x,y) &=&\frac{\sum_{i=-\infty }^{\infty }\sum_{j=-\infty }^{\infty }\Psi
(i,j)x^{i}y^{j}}{1-\sum_{u=1}^{n}\sum_{v=1}^{m}c_{uv}x^{s-(u-1)}y^{t-(v-1)}z}
\\
&=&\sum_{i=-\infty }^{\infty }\sum_{j=-\infty }^{\infty }\Psi
(i,j)x^{i}y^{j}\sum_{k=0}^{\infty
}(\sum_{u=1}^{n}\sum_{v=1}^{m}c_{uv}x^{s-(u-1)}y^{t-(v-1)})^{k}z^{k} \\
&=&\sum_{k=0}^{\infty }\sum_{i=-\infty }^{\infty }\sum_{j=-\infty }^{\infty
}\Psi
(i,j)x^{i}y^{j}(%
\sum_{u=1}^{n}c_{u1}x^{s-(u-1)}y^{t}+c_{u2}x^{s-(u-1)}y^{t-1}+...+ \\
&&c_{um}x^{s-(u-1)}y^{t-(m-1)})^{k}z^{k} \\
&=&\sum_{k=0}^{\infty }\sum_{i=-\infty }^{\infty }\sum_{j=-\infty }^{\infty
}\Psi
(i,j)x^{i}y^{j}(c_{11}x^{s}y^{t}+c_{12}x^{s}y^{t-1}+...+c_{1m}x^{s}y^{t-(m-1)}+
\\
&&c_{21}x^{s-1}y^{t}+c_{22}x^{s-1}y^{t-1}+...+c_{2m}x^{s-1}y^{t-(m-1)}+...+
\\
&&c_{n1}x^{s-(n-1)}y^{t}+c_{n2}x^{s-(n-1)}y^{t-1}+...+c_{nm}x^{s-(n-1)}y^{t-(m-1)})^{k}z^{k}
\\
&=&\sum_{k=0}^{\infty }\sum_{i=-\infty }^{\infty }\sum_{j=-\infty }^{\infty
}\Psi
(i,j)x^{i}y^{j}c_{11}^{s_{11}}x^{ss_{11}}y^{ts_{11}}c_{12}^{s_{12}}x^{ss_{2}}y^{(t-1)s_{12}}...c_{1m}^{s_{m}}x^{ss_{1m}}y^{[t-(m-1)]s_{1m}}
\\
&&c_{21}^{s_{21}}x^{(s-1)s_{21}}y^{ts_{21}}c_{22}^{s_{22}}x^{(s-1)s_{22}}y^{(t-1)s_{22}}...c_{2m}^{s_{2m}}x^{(s-1)s_{2m}}y^{[t-(m-1)]s_{2m}}...
\\
&&c_{n1}^{s_{n1}}x^{[s-(n-1)]s_{n1}}y^{ts_{n1}}c_{n2}^{s_{n2}}x^{[s-(n-1)]s_{n2}}y^{(t-1)s_{n2}}...c_{nm}^{s_{nm}}x^{[s-(n-1)]s_{nm}}y^{[t-(m-1)]s_{nm}}z^{k}
\\
&=&\sum_{k=0}^{\infty }\sum_{i=-\infty }^{\infty }\sum_{j=-\infty }^{\infty
}\Psi (i,j)x^{i}y^{j}\sum_{s_{qr\in I}}\binom{k}{%
s_{11},s_{12},...,s_{1m},...,s_{nm}}\prod_{\substack{u=1,2,...,n,\\v=1,2,...,m}} 
c_{uv}^{s_{uv}}x^{a}y^{b}z^{k} \\
&=&\sum_{k=0}^{\infty }\sum_{i=-\infty }^{\infty }\sum_{j=-\infty }^{\infty
}\sum_{s_{r\in I}}\binom{k}{s_{11},s_{12},...,s_{1m},...,s_{nm}}\Psi
(i,j)\prod_{\substack{u=1,2,...,n,\\v=1,2,...,m}}%
c_{uv}^{s_{uv}}x^{a+i}y^{b+j}z^{k},
\end{eqnarray*}
where 
\begin{eqnarray*}
I &=&\{s_{11},s_{12},...,s_{1m},...,s_{nm}|0\leq s_{uv}\leq
k,\sum_{u=1}^{n}\sum_{v=1}^{m}s_{uv}=k\}, \\
a
&=&s\sum_{v=1}^{m}s_{1v}+(s-1)\sum_{v=1}^{m}s_{2v}+...+[s-(n-1)]%
\sum_{v=1}^{m}s_{nv} \\
&=&\sum_{u=1}^{n}[s-(u-1)]\sum_{v=1}^{m}s_{uv} \\
&=&\sum_{u=1}^{n}\sum_{v=1}^{m}[s-(u-1)]s_{uv} \\
&=&\sum_{u=1}^{n}\sum_{v=1}^{m}[(s+1)-u)]s_{uv} \\
&=&(s+1)k-\sum_{u=1}^{n}\sum_{v=1}^{m}us_{uv}\text{, } \\
b &=&\sum_{u=1}^{n}\sum_{v=1}^{m}[t-(v-1)]s_{uv} \\
&=&\sum_{u=1}^{n}\sum_{v=1}^{m}[(t+1)-v]s_{uv} \\
&=&(t+1)k-\sum_{u=1}^{n}\sum_{v=1}^{m}vs_{uv}.
\end{eqnarray*}
Let 
\[
a+i=l\text{, }b+j=p. 
\]
Then we have

\[
h(x,y,z)=\sum_{k=0}^{\infty }\sum_{l=-\infty }^{\infty }\sum_{p=-\infty
}^{\infty }\sum_{s_{qr\in I}}\binom{k}{s_{11},s_{12},...,s_{1m},...,s_{nm}}%
\Psi (l-a,p-b)\prod_{\substack{u=1,2,...,n,\\v=1,2,...,m}}%
c_{uv}^{s_{uv}}x^{l}y^{p}z^{k}. 
\]
So the solution is of the form 
\[
U_{i,j}=\sum_{s_{qr\in I}}\binom{k}{s_{11},s_{12},...s_{1m},...,s_{nm}}\Psi
(i-a,j-b)\prod_{\substack{u=1,2,...,n,\\v=1,2,...,m}}c_{uv}^{s_{uv}}, 
\]
\endproof%

Now, we extend eq(\ref{eq7}) to n-dimensional case, the high-order partial difference equation (HOPDE), as follows, and find out it's exact solution.
\begin{lemma}
\begin{flalign*}
(\sum_{u_{1}=1}^{k_{1}}\sum_{u_{2}=1}^{k_{2}}...%
\sum_{u_{n-1}=1}^{k_{n-1}}
c_{u_{1}u_{2}...u_{n-1}}%
\prod_{i=1}^{n-1}x_{i}^{s_{i}-(u_{i}-1)})^{j}
=\sum_{r_{u_{1}u_{2}...u_{n-1}}\in I}\binom{j}{%
r_{11...1},r_{11...2},...,r_{11...k_{n-1}},...,r_{k_{1}k_{2}...k_{n-1}}}\\
\prod_{\substack{u_{1}=1,2,...,k_{1}.\\u_{2}=1,2,...,k_{2.}\\...\\u_{n-1}=1,2,...,k_{n-1.}}}
c_{u_{1}u_{2}...u_{n-1}}^{r_{u_{1}u_{2}...u_{n-1}}}%
\prod_{i=1}^{n-1}x_{i}^{a_{i}},
\end{flalign*}
where

\begin{flalign*}
I={r_{11...1},r_{11...2},...,r_{11...k_{n-1}},...,r_{k_{1}k_{2}...k_{n-1}}|%
0 \leq r_{u_{1}u_{2}...u_{n-1}}}\leq
j,\sum_{u_{1}=1}^{k_{1}}\sum_{u_{2}=1}^{k_{2}}...%
\sum_{u_{n-1}}^{k_{n-1}}r_{u_{1}u_{2}...u_{n-1}}=j\}, \\
a_{i}=(s_{i}+1)j-\sum_{u_{1}=1}^{k_{1}}\sum_{u_{2}=1}^{k_{2}}...%
\sum_{u_{n-1}=1}^{k_{n-1}}u_{i}r_{u_{1}u_{2}...u_{n-1}},\forall i=1,2,...n-1.
\end{flalign*}
\end{lemma}
\proof%
\begin{flalign*}
(\sum_{u_{1}=1}^{k_{1}}\sum_{u_{2}=1}^{k_{2}}...%
\sum_{u_{n-1}=1}^{k_{n-1}}c_{u_{1}u_{2}...u_{n-1}}%
\prod_{i=1}^{n-1}x_{i}^{s_{i}-(u_{i}-1)})^{j} \\
=\sum_{u_{1}=1}^{k_{1}}\sum_{u_{2}=1}^{k_{2}}...%
\sum_{u_{n-2}=1}^{k_{n-2}}(c_{u_{1}u_{2}...u_{n-2}1}%
\prod_{i=1}^{n-2}x_{i}^{s_{i}-(u_{i}-1)}x_{n-1}^{s_{n}}+c_{u_{1}u_{2}...u_{n-2}2}\prod_{i=1}^{n-2}x_{i}^{s_{i}-(u_{i}-1)}x_{n-1}^{s_{n}-1}+
\\
...+c_{u_{1}u_{2}...u_{n-2}k_{n-1}}%
\prod_{i=1}^{n-2}x_{i}^{s_{i}-(u_{i}-1)}x_{n-1}^{s_{n}-(k_{n}-1)})]^{j} \\
=[c_{11...11}\prod_{i=1}^{n-1}x_{i}^{s_{i}}+c_{11...12}%
\prod_{i=1}^{n-2}x_{i}^{s_{i}-1}x_{n-1}^{s_{n-1}-1}+...+c_{11...1k_{n-1}}%
\prod_{i=1}^{n-2}x_{i}^{s_{i}-1}x_{n-1}^{s_{n-1}-(k_{n-1}-1)}+ \\
...+c_{k_{1}k_{2}...k_{n-1}}\prod_{i=1}^{n-1}x_{i}^{s_{i}-(k_{i}-1)}]^{j}
\\
=\sum_{r_{u_{1}u_{2}...u_{n-1}}\in I}\binom{j}{%
r_{11...1},r_{11...12},...,r_{11...1k_{n-1}},...,r_{k_{1}k_{2}...k_{n-1}}}%
(c_{11...11}\prod_{i=1}^{n-1}x_{i}^{s_{i}})^{r_{11...11}} \\
(c_{11...12}%
\prod_{i=1}^{n-2}x_{i}^{s_{i}-1}x_{n-1}^{s_{n-1}-1})^{r_{11...12}}(c_{11...1k_{n-1}}\prod_{i=1}^{n-2}x_{i}^{s_{i}-1}x_{n-1}^{s_{n-1}-(k_{n-1}-1)})^{r_{11...1k_{n-1}}}
\\
...(c_{k_{1}k_{2}...k_{n-1}}%
\prod_{i=1}^{n-1}x_{i}^{s_{i}-(k_{i}-1)})^{r_{k_{1}k_{2}...k_{n-1}}} \\
=\sum_{r_{u_{1}u_{2}...u_{n-1}}\in I}\binom{j}{%
r_{11...1},r_{11...2},...,r_{11...k_{n-1}},...,r_{k_{1}k_{2}...k_{n-1}}} %
\prod_{\substack{u_{1}=1,2,...,k_{1}.\\u_{2}=1,2,...,k_{2.}\\...\\u_{n-1}=1,2,...,k_{n-1.}}}\\%
c_{u_{1}u_{2}...u_{n-1}}^{r_{u_{1}u_{2}...u}{}_{n-1}}%
\prod_{i=1}^{n-1}x_{i}^{a_{i}},
\end{flalign*}
where

\begin{eqnarray*}
I
&=&%
\{r_{11...1},r_{11...12},...,r_{11...1k_{n-1}},...,r_{k_{1}k_{2}...k_{n-1}}|
\\
0 &\leq &r_{u_{1}u_{2}...u_{n-1}}\leq
j,\sum_{u_{1}=1}^{k_{1}}\sum_{u_{2}=1}^{k_{2}}...%
\sum_{u_{n-1}}^{k_{n-1}}r_{u_{1}u_{2}...u_{n-1}}=j\}, \\
\forall i &=&1,2,...,n-1, \\
a_{i}
&=&\sum_{u_{1}=1}^{k_{1}}\sum_{u_{2}=1}^{k_{2}}...%
\sum_{u_{n-1}=1}^{k_{n-1}}[s_{i}-(u_{i}-1)]r_{u_{1}u_{2}...u_{n-1}} \\
&=&\sum_{u_{1}=1}^{k_{1}}\sum_{u_{2}=1}^{k_{2}}...%
\sum_{u_{n-1}=1}^{k_{n-1}}[(s_{i}+1)-u_{i}]r_{u_{1}u_{2}...u_{n-1}} \\
&=&(s_{i}+1)j-\sum_{u_{1}=1}^{k_{1}}\sum_{u_{2}=1}^{k_{2}}...%
\sum_{u_{n-1}=1}^{k_{n-1}}u_{i}r_{u_{1}u_{2}...u_{n-1}}.
\end{eqnarray*}
\endproof%
\begin{theorem}
Let 
\[
U_{e_{1}+s_{1},e_{2}+s_{2},...,e_{n-1}+s_{n-1},e_{n}+1}=%
\sum_{u_{1}=1}^{k_{1}}\sum_{u_{2}=1}^{k_{2}}...%
\sum_{u_{n-1}=1}^{k_{n-1}}c_{u_{1}u_{2}...u_{n-1}}U_{e_{1}+u_{1}-1,e_{2}+u_{2}-1,...,e_{n-1}+u_{n-1}-1,e_{n}} 
\]
be a linear difference equation in e$_{1}$, e$_{2}$,...,e$_{n}$, where e$%
_{i}\in Z$ for i$\in \{1,2,...,n-1\}$, $e_{n}\in N\cup \{0\}$. \ Given $%
U_{e_{1},e_{2},...,e_{n-1},0}=\Psi (e_{1},e_{2},...,e_{n-1})$, then its
solution is of the form 
\begin{eqnarray*}
U_{e_{1},e_{2},...,e_{n}} &=&\sum_{r_{u_{1}u_{2}...u_{n-1}}\in I}\binom{j}{%
r_{11...1},r_{11...2},...,r_{11...k_{n-1}},...,r_{k_{1}k_{2}...k_{n-1}}} \\
&&\Psi (e_{1}-a_{1},e_{2}-a_{2},...,e_{n-1}-a_{n-1})\prod_{\substack{u_{1}=1,2,...,k_{1}.\\u_{2}=1,2,...,k_{2.}\\ ...\\u_{n-1}=1,2,...,k_{n-1.}}}%
c_{u_{1}u_{2}...u_{n-1}}^{r_{u_{1}u_{2}...u}{}_{n-1}},
\end{eqnarray*}
\end{theorem}

where 
\begin{eqnarray*}
I
&=&\{r_{11...1},r_{11...2},...,r_{11...k_{n-1}},...,r_{k_{1}k_{2}...k_{n-1}}|
\\
0 &\leq &r_{u_{1}u_{2}...u_{n-1}}\leq
j,\sum_{u_{1}=1}^{k_{1}}\sum_{u_{2}=1}^{k_{2}}...%
\sum_{u_{n-1}}^{k_{n-1}}r_{u_{1}u_{2}...u_{n-1}}=j\}, \\
a_{i}
&=&(s_{i}+1)j-\sum_{u_{1}=1}^{k_{1}}\sum_{u_{2}=1}^{k_{2}}...%
\sum_{u_{n-1}=1}^{k_{n-1}}u_{i}r_{u_{1}u_{2}...u_{n-1}},\forall i=1,2,...n-1.
\end{eqnarray*}

\proof%
Let 
\[
\stackrel{\rightharpoonup }{x}=(x_{1},x_{2},...,x_{n}), 
\]
\[
h(\stackrel{\rightharpoonup }{x})=\sum_{e_{1}=-\infty }^{\infty
}\sum_{e_{2}=-\infty }^{\infty }...\sum_{e_{n-1}=-\infty }^{\infty
}\sum_{e_{n}=0}^{\infty
}U_{e_{1},e_{2},...,e_{n}}\prod_{i=1}^{n}x_{i}^{e_{i}}. 
\]
Then we have

\begin{eqnarray*}
&&\sum_{e_{1}=-\infty }^{\infty }\sum_{e_{2}=-\infty }^{\infty
}...\sum_{e_{n-1}=-\infty }^{\infty }\sum_{e_{n}=0}^{\infty
}U_{e_{1}+s_{1},e_{2}+s_{2},...,e_{n-1}+s_{n-1},e_{n}+1}%
\prod_{i=1}^{n}x_{i}^{e_{i}} \\
&=&\sum_{e_{1}=-\infty }^{\infty }\sum_{e_{2}=-\infty }^{\infty
}...\sum_{e_{n-1}=-\infty }^{\infty }\sum_{e_{n}=0}^{\infty
}\sum_{u_{1}=1}^{k_{1}}\sum_{u_{2}=1}^{k_{2}}...\sum_{u_{n-1}=1}^{k_{n-1}} \\
&&c_{u_{1}u_{2}...u_{n-1}}U_{e_{1}+u_{1}-1,e_{2}+u_{2}-1,...,e_{n-1}+u_{n-1}-1,e_{n}}\prod_{i=1}^{n}x_{i}^{e_{i}},
\end{eqnarray*}

\begin{eqnarray*}
&&\frac{1}{\prod_{i=1}^{n-1}x_{i}^{s_{i}}x_{n}}[h(\stackrel{\rightharpoonup 
}{x})-\sum_{e_{1}=-\infty }^{\infty }\sum_{e_{2}=-\infty }^{\infty
}...\sum_{e_{n-1}=-\infty }^{\infty }\Psi
(e_{1},e_{2},...,e_{n-1})\prod_{i=1}^{n-1}x_{i}^{e_{i}}] \\
&=&\sum_{u_{1}=1}^{k_{1}}\sum_{u_{2}=1}^{k_{2}}...%
\sum_{u_{n-1}=1}^{k_{n-1}}c_{u_{1}u_{2}...u_{n-1}}%
\prod_{i=1}^{n-1}x_{i}^{-(u_{i}-1)}h(\stackrel{\rightharpoonup }{x}),
\end{eqnarray*}

\begin{eqnarray*}
&&h(\stackrel{\rightharpoonup }{x})(1-\sum_{u_{1}=1}^{k_{1}}%
\sum_{u_{2}=1}^{k_{2}}...\sum_{u_{n-1}=1}^{k_{n-1}}c_{u_{1}u_{2}...u_{n-1}}%
\prod_{i=1}^{n-1}x_{i}^{s_{i}-(u_{i}-1)}x_{n}) \\
&=&\sum_{e_{1}=-\infty }^{\infty }\sum_{e_{2}=-\infty }^{\infty
}...\sum_{e_{n-1}=-\infty }^{\infty }\Psi
(e_{1},e_{2},...,e_{n-1})\prod_{i=1}^{n-1}x_{i}^{e_{i}},
\end{eqnarray*}

\begin{eqnarray*}
h(\stackrel{\rightharpoonup }{x}) &=&\frac{\sum_{e_{1}=-\infty }^{\infty
}\sum_{e_{2}=-\infty }^{\infty }...\sum_{e_{n-1}=-\infty }^{\infty }\Psi
(e_{1},e_{2},...,e_{n-1})\prod_{i=1}^{n-1}x_{i}^{e_{i}}}{(1-%
\sum_{u_{1}=1}^{k_{1}}\sum_{u_{2}=1}^{k_{2}}...%
\sum_{u_{n-1}=1}^{k_{n-1}}c_{u_{1}u_{2}...u_{n-1}}%
\prod_{i=1}^{n-1}x_{i}^{s_{i}-(u_{i}-1)}x_{n})} \\
&=&\sum_{e_{1}=-\infty }^{\infty }\sum_{e_{2}=-\infty }^{\infty
}...\sum_{e_{n-1}=-\infty }^{\infty }\Psi
(e_{1},e_{2},...,e_{n-1})\prod_{i=1}^{n-1}x_{i}^{e_{i}} \\
&&\sum_{j=0}^{\infty
}(\sum_{u_{1}=1}^{k_{1}}\sum_{u_{2}=1}^{k_{2}}...%
\sum_{u_{n-1}=1}^{k_{n-1}}c_{u_{1}u_{2}...u_{n-1}}%
\prod_{i=1}^{n-1}x_{i}^{s_{i}-(u_{i}-1)})^{j}x_{n}^{j}.
\end{eqnarray*}
Applying above lemma, we have

\begin{eqnarray*}
h(\stackrel{\rightharpoonup }{x}) &=&\sum_{j=0}^{\infty }\sum_{e_{1}=-\infty
}^{\infty }\sum_{e_{2}=-\infty }^{\infty }...\sum_{e_{n-1}=-\infty }^{\infty
}\Psi (e_{1},e_{2},...,e_{n-1})\prod_{i=1}^{n-1}x_{i}^{e_{i}} \\
&&\sum_{r_{u_{1}u_{2}...u_{n-1}}\in I}\binom{j}{%
r_{11...1},r_{11...2},...,r_{11...k_{n-1}},...,r_{k_{1}k_{2}...k_{n-1}}} \\
&&\prod_{\substack{u_{1}=1,2,...,k_{1}.\\u_{2}=1,2,...,k_{2.}\\...\\u_{n-1}=1,2,...,k_{n-1.}}}%
c_{u_{1}u_{2}...u_{n-1}}^{r_{u_{1}u_{2}...u}{}_{n-1}}%
\prod_{i=1}^{n-1}x_{i}^{a_{i}}x_{n}^{j} \\
&=&\sum_{j=0}^{\infty }\sum_{e_{1}=-\infty }^{\infty }\sum_{e_{2}=-\infty
}^{\infty }...\sum_{e_{n-1}=-\infty }^{\infty
}\sum_{r_{u_{1}u_{2}...u_{n-1}}\in I}\binom{j}{%
r_{11...1},r_{11...2},...,r_{11...k_{n-1}},...,r_{k_{1}k_{2}...k_{n-1}}} \\
&&\Psi (e_{1},e_{2},...,e_{n-1})\prod_{\substack{
u_{1}=1,2,...,k_{1}.\\u_{2}=1,2,...,k_{2.}\\...\\u_{n-1}=1,2,...,k_{n-1.}}}%
c_{u_{1}u_{2}...u_{n-1}}^{r_{u_{1}u_{2}...u}{}_{n-1}}%
\prod_{i=1}^{n-1}x_{i}^{a_{i}+e_{i}}x_{n}^{j}.
\end{eqnarray*}
Let 
\[
a_{i}+e_{i}=b_{i},\forall i=1,2,...,n-1. 
\]
Then we have

\begin{eqnarray*}
h(\stackrel{\rightharpoonup }{x})=\sum_{j=0}^{\infty }\sum_{b_{1}=-\infty
}^{\infty }\sum_{b_{2}=-\infty }^{\infty }...\sum_{b_{n-1}=-\infty }^{\infty
}\sum_{r_{u_{1}u_{2}...u_{n-1}}\in I}\binom{j}{%
r_{11...1},r_{11...2},...,r_{11...k_{n-1}},...,r_{k_{1}k_{2}...k_{n-1}}} \\
\Psi (b_{1}-a_{1},b_{2}-a_{2},...,b_{n-1}-a_{n-1})\prod_{\substack{%
u_{1}=1,2,...,k_{1}.\\u_{2}=1,2,...,k_{2.}\\...\\u_{n-1}=1,2,...,k_{n-1.}}}
c_{u_{1}u_{2}...u_{n-1}}^{r_{u_{1}u_{2}...u}{}_{n-1}}%
\prod_{i=1}^{n-1}x_{i}^{b_{i}}x_{n}^{j},
\end{eqnarray*}
Thus the solution is of the form 
\begin{eqnarray*}
U_{e_{1},e_{2},...,e_{n}}=\sum_{r_{u_{1}u_{2}...u_{n-1}}\in I}\binom{j}{%
r_{11...1},r_{11...2},...,r_{11...k_{n-1}},...,r_{k_{1}k_{2}...k_{n-1}}} \\
\Psi (e_{1}-a_{1},e_{2}-a_{2},...,e_{n-1}-a_{n-1})\prod_{\substack{%
u_{1}=1,2,...,k_{1}.\\u_{2}=1,2,...,k_{2.}\\...\\u_{n-1}=1,2,...,k_{n-1.}}}%
c_{u_{1}u_{2}...u_{n-1}}^{r_{u_{1}u_{2}...u}{}_{n-1}},
\end{eqnarray*}
\endproof%

\bibliographystyle{amsplain}
\input pde.new.ref.bbl

\end{document}

%% file: pde.new.ref.bbl
\providecommand{\bysame}{\leavevmode\hbox to3em{\hrulefill}\thinspace}
\providecommand{\MR}{\relax\ifhmode\unskip\space\fi MR }
\providecommand{\MRhref}[2]{%
  \href{http://www.ams.org/mathscinet-getitem?mr=#1}{#2}
}
\providecommand{\href}[2]{#2}